\documentclass{article}
\usepackage{amssymb,amsmath,amsfonts,enumerate,url,graphicx}

\newtheorem{theorem}{Theorem}[section]

\newtheorem{definition}[theorem]{Definition}

\begin{document}

\title{Pedal curves of conics: an automated exploration of some cubics, sextics, octics and more}

\author{Thierry Dana-Picard}

\maketitle

\begin{center}
Department  of  Mathematics,  Jerusalem  College  of
Technology\\  
Jerusalem 9116011, Israel\\e-mail: ndp@jct.ac.il
\end{center}



\begin{abstract}
Constructions and exploration of plane algebraic curves has received a new push with the development of automated methods, whose algorithms are continuously improved and implemented in various software packages. We use them to explore the pedal curves of conics. This provides a construction of interesting geometric loci, given at first by parametric representations. After implicitization, they appear as sextics, octics and other curves of higher degree. We explore their irreducibility and their singular points (crunodes, cusps, etc.).  
\end{abstract}

\maketitle
\date{\today}

\section{Introduction}
Plane algebraic curves are a traditional topic, which drew continuous interest, from ancient times to nowadays. A vast knowledge is offered in books such as \cite{salmon,yates} and websites \cite{mathcurve}.Algorithms to work in polynomial rings (computation of Gr\" obner bases, elimination, etc.) led to both the revival of traditional topics  and to the exploration of varieties in @d and 3D spaces (and beyond) which have been somehow forgotten or at least left aside. Among these are envelopes of families of curves dependent on a parameter \cite{revival,envelopes and offsets}, isoptic and bisoptic curves \cite{DPK-2018,DMZ,isoptics of astroid,DZM,DNMC,inner} and also constructions and exploration of curves of higher degree.  Some of them appear as an outcome of a "pure geometric" question, we mean geometric loci \cite{octics,thales etc.,thales etc. MT}. For example, we can mention hyperbolisms of curves with respect to a point and a line \cite{hyperbolisms}, which have been rarely studied. 

The exploration is  based on the usage of software, originally called either Computer Algebra Systems (CAS) or Dynamic Geometry Software (DGS). Different software may have similar features, and the distinction between CAS and DGS became quite fuzzy, as they evolved to multi-purpose software. Anyway, the different possibilities offered by the different packages make a dialog between them very useful. Networking between them is of the utmost importance \cite{dialog}. A pitfall is that transferring data from one pckage to another one must often be made "by hand", i.e. by copy-paste using the mouse. An automatic dialog between them is a goal mentioned a long time ago \cite{eugenio}, but has not been achieved yet.

In the present work, we explore \emph{pedal curves} of conics, namely parabolas and ellipses. As the reader will see later, these curves can be realized as geometric loci (according to Definition \ref{def pedal}, but also as envelopes of some families of circles (see Section \ref{pedal of parabola}).

\begin{definition}
\label {def pedal}
Let be given a plane curve $\mathcal{C}$ and a point $D$.  The \emph{pedal} of $\mathcal{C}$  with respect to $D$ (called the \emph{pole}) is the geometric locus of the feet of the lines passing by $D$ and perpendicular to the tangents to $\mathcal{C}$.
\end{definition}
We will generally denote that curve by $\mathcal{P}$.

In this work, we observe and explore curves of degree 3 (cubics), 6 (sextics), 8 (octics) and in some cases curves of degree 20. These will be reducible. A central issue is to factor polynomials in order to check whether the obtained pedal curves are irreducible or not, and what happens if superfluous components appear. Such an issue has already been discussed in \cite{isoptics of quartics}.

We look at the curve under study as an intriguing curve, in the spirit of \cite{ferrarello et al}. The construction process provides a parametric representation, which is implicitized to obtain a polynomial equation. Different methods exist, based on resultants or on Gr\"obner bases \cite{cox}, or more recently on the Gr\"obner cover \cite{montes}. In particular, this book displays numerous examples of determination of geometric loci, which is our concern here. Elimination is almost ubiquitous in the papers mentioned above, and algorithms for elimination are implemented in both kinds of software that we use: 
\begin{itemize}
\item  The Computer Algebra System Maple 2024, in particular for its strong \emph{PolynomialIdeals} package. Maple has also an \emph{algcurves} packaage, with a command for automatic implicitization, but we preferred not to use it,as it is valid in specific cases.
\item GeoGebra Discovery (GD) \footnote{Freely downloadable from \url{https://github.com/kovzol/geogebra-discovery}. We recommend to check frequently for the newest release.}, a package for automated commands, which works on the basis of the regular GeoGebra \footnote{Freely downloadable from \url{http://www.geogebra.org}}. 
\end{itemize}
GD provides accurate answers for the determination of a geometric locus (in several settings), based on symbolic algorithms. It may provide non only a plot, but under certain conditions, also a polynomial equation for the geometric locus. An embedded Computer Algebra System called Giac \cite{kovacs parisse} is helpful to decide whether the determined locus is reducible or not, but sometimes it is not strong enough, whence the need to use Maple. This will be illustrated in Section \ref{section ellipses}.

\section{Pedal curves of a parabola}
\label{pedal of parabola}

In what follows, we work out specific example, but they are good representatives of a general situation. We chose the examples for the output to be not too complicated, and the topoology of the curves is clear enough on the plots. Of course, proofs are needed for what can be conjectured form the visualization.

\subsection{Exploration with GeoGebra Discovery}
WLOG, we consider a parabola $\mathcal{C}$, with focus $F(0,1)$ and whose directrix has equation $y=-1$. This parabola has equation $x^2-4y=0$.
Many commands require the construction to be purely geometric. Therefore, we make the following construction (the notations appear in the figures):
\begin{enumerate}[(i)]
\item Plot two points $A$ and $B$ to determine a line, which will be the directrix of the parabola, denoted by $c$.
\item Plot a point $F$ not on the line $AB$, the focus of the parabola.
\item The command \textbf{Locus($<$Point defining the locus$>$, $<$Moving Point$>$)} plots the parabola.
\item Plot a point $E$  on the parabola, with the button \emph{Point on Object} (a written command is also available).
\item Determine the tangent to the parabola at $E$.
\item Plot any point $D$ in the plane.
\item Determine the normal $h$ through $D$ to the tangent $g$ to the parabola at $E$.
\item determine the point of intersection $H$ of this normal with $g$.
\item The command \textbf{Locus(\emph{H,E})} provides a plot of the pedal curve of the parabola with respect to the point $E$.
\item The command \textbf{LocusEquation(\emph{H,E})} provides a symbolic equation for this curve.
\end{enumerate}

\begin{figure}[htb]
\begin{center}
\includegraphics[width=3.5cm]{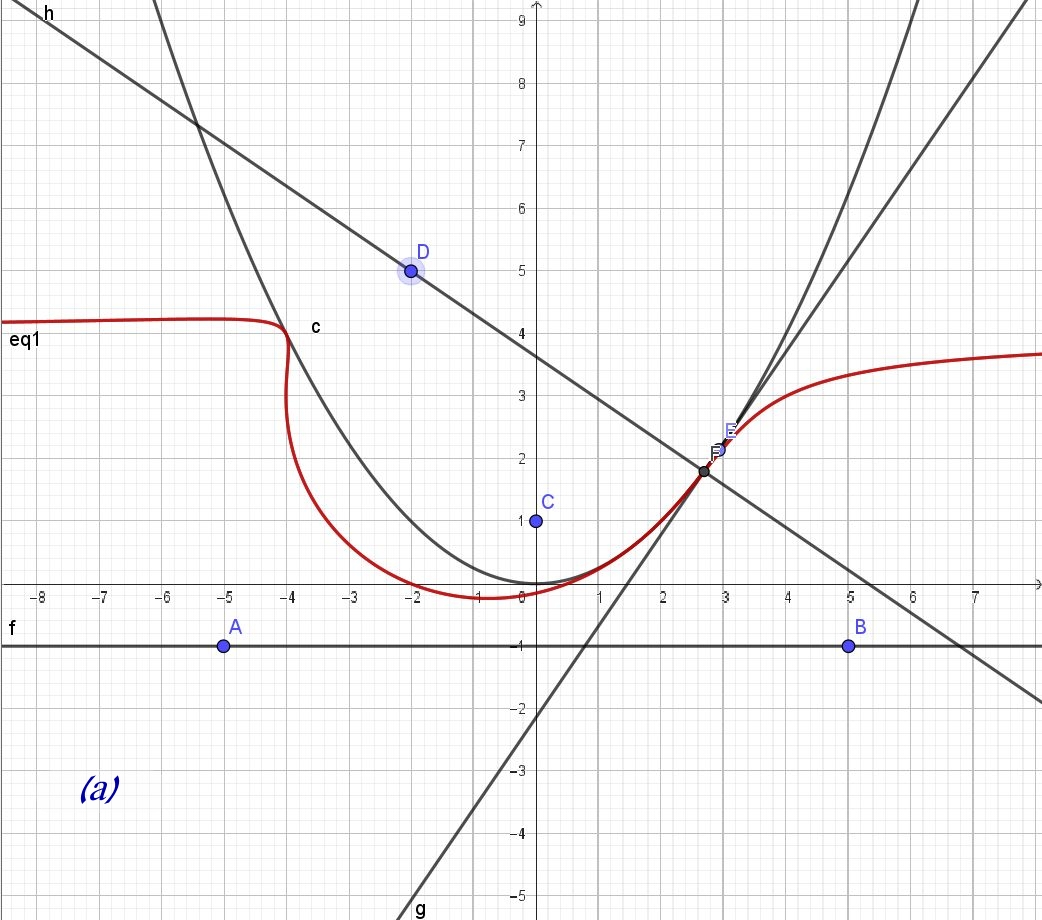}
\quad
\includegraphics[width=3.5cm]{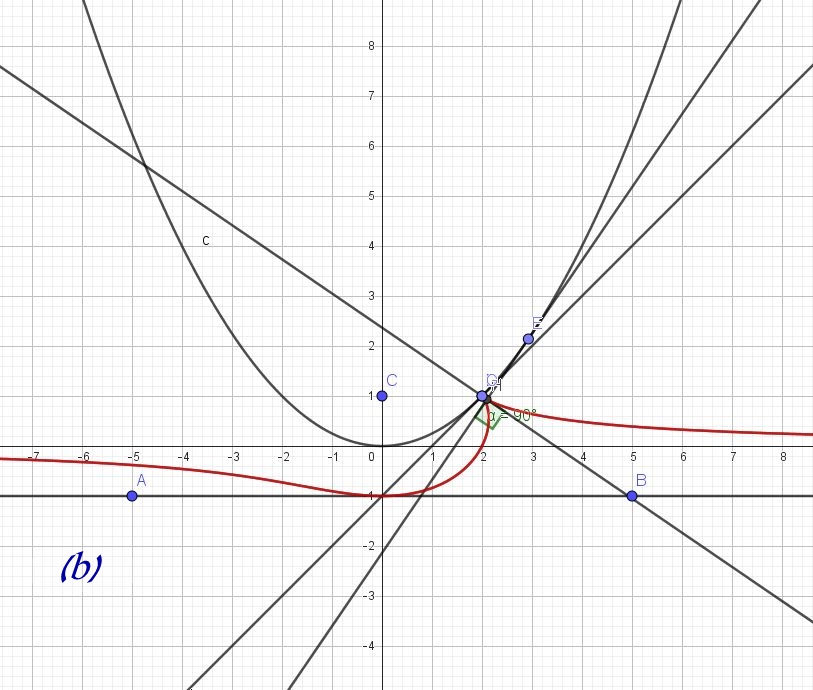}
\quad
\includegraphics[width=3.5cm]{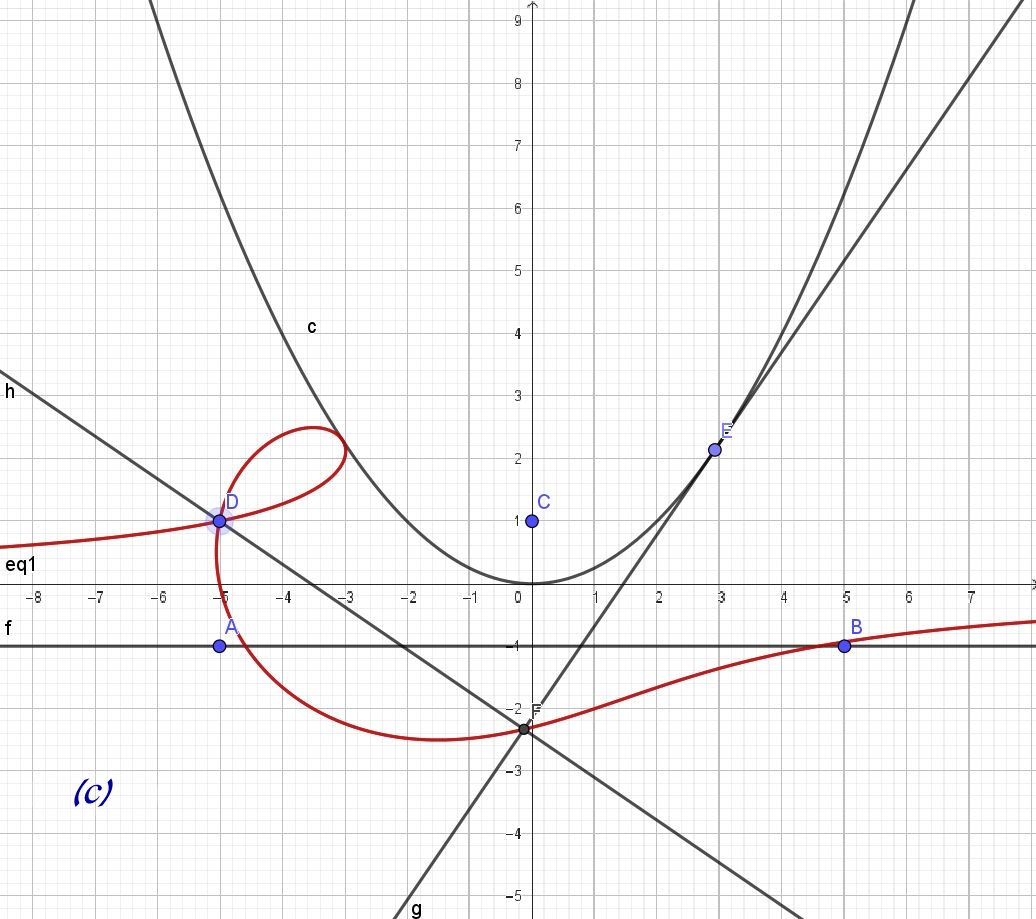}
\caption{Pedal curve of a canonical parabola with respect to a point - the 3 cases}
\label{fig pedal of a parabola}
\end{center}
\end{figure}

Figure \ref{fig pedal of a parabola} shows screenshots of the exploration.
The differential properties of the obtained curves can be proven using symbolic equations, which are derived from the algebraic work in next subsection. 
The obtained curves show 3 different topologies:
\begin{enumerate}[(i)]
\item If $D$ is a point inside the parabola, the pedal curve is a smooth cubic, as shown in Figure \ref{fig pedal of a parabola}(a). 
\item Id $D$ is a point on the parabola, then the pedal curve has a cusp at $D$; see  Figure \ref{fig pedal of a parabola}(b)
\item   If $D$ is a point out of the parabola, the pedal curve is a crunodal cubic, i.e. it has a crunode (a double point, where the curve intersects itself with two different tangents; see \cite{salmon}), as shown in Figure \ref{fig pedal of a parabola}(c).
\end{enumerate}
These properties can be proven in general. In order not to display too heavy expressions, we study examples in the next section. We wish to mention that the polynomial equations in these examples are obtained via algebraic computations with Maple, but they have also been obtained with GeoGebra's \textbf{LocusEquation} command.
 
\subsection{Algebraic work with Maple}
\label{subsection alg work for parabola}
WLOG, we consider a parabola $\mathcal{C}$, whose focus is the point $F(0,1)$ and the directrix has equation $y=-1$. This parabola has equation $x^2-4y=0$. At a point $A_0=(x_0,y_0)$ on the parabola, the tangent $T_{A_0}$ has equation $y = \frac{1}{2}xx_0 - \frac 14 x_0^2$. Denote $D=(x_D,y_D)$; then the perpendicular to $T_{A_0}$ through $D$ has equation $y - y_D = -\frac{2}{x_0}(x - x_D)$. The foot $H$ of this perpendicular on the tangent $T_{A_0}$ is thus the point verifying the following system of equations:
\begin{equation}
\label{eq foot of the perpendicular}
\begin{cases}
y = \frac{1}{2}xx_0 - \frac 14 x_0^2\\
y - y_D = -\frac{2}{x_0}(x - x_D)
\end{cases}
\end{equation}
The solution is given by:
\begin{equation}
\label{coord pedal point 1}
\begin{cases}
x_H = \frac{x_0^3 + 4x_0y_D + 8x_D}{2(x_0^2 + 4)}\\
y_H = \frac{x_0(x_0y_D - x_0 + 2x_D)}{x_0^2 + 4}
\end{cases}.
\end{equation}
Maybe the general formula is not clear enough to be an incitement to distinguish three cases, but the dynamic experiment with the DGS provides it.

There exist infinitely many parametric representations for the given parabola; we take $(x,y)=(2t,t^2),\; t \in \mathbb{R}$. Thus we have:
\begin{equation}
\label{coord pedal point 2}
\begin{cases}
x_H =  \frac{8t^3 + 8ty_D + 8x_D}{2(4t^2 + 4)}\\
y_H = \frac{2t(2ty_D - 2t + 2x_D)}{4t^2 + 4}
\end{cases}.
\end{equation}
From this point, we work with Maple (version 2024). Define two polynomials and the ideal that they generate, using the commands: 
\small
\begin{verbatim}
parab:= implicitplot(4y=x^2,x=-8..8,y=-3..8,color=blue);
tgt := y = 1/2*x*x0 - 1/4*x0^2;
nrml := y - yD = -2*(x - xD)/x0;
H := solve({nrml, tgt}, {x, y})
H := subs(x0 = 2*t, subs(y0 = t^2, H));
P1 := x*denom(rhs(H[1])) - numer(rhs(H[1])):
P2 := y*denom(rhs(H[2])) - numer(rhs(H[2]));
J:=<P1,P2>;
\end{verbatim}
\normalsize
We eliminate the parameter $t$ using the \emph{PolynomialIdeal} and discover that the obtained polynomial is generated by a unique polynomial of degree 3. Then a single command shows that the generator is irreducible. The code is as follows: 
\begin{verbatim}
JE := EliminationIdeal(J, {x, xD, y, yD});
ped := Generators(JE)[1];
factors(ped);
\end{verbatim}
The generator of the ideal (denoted by $ped$ in the code) provides a symbolic equation for the desired pedal curve is:
\footnotesize
\begin{equation}
\label{eq pedal of the parabola}
\begin{aligned}
G(x,y) & = (x_D^2+y_D^2-2y_D+1) \; y^3+ (x_D^2+y_D^2-2y_D+1) x^2y    \\
\quad &  - (x_D^2 y_D+y_D^3+x_D^2-3y_D^2+3y_D+1)\; x^2 - (x_D^3 + x_D  y_D^2-2x_Dy_D+x_D) \; xy     \\
\quad &   - (2 x_D^2 y_D+2y_D^3-4y_D^2+2y_D) \; y^2  + (x_D^2 y_D^2+y_D^4-2y_D^3+y_D^2) \; y\\
\quad & + (x_D^3 y_D+x_D y_D^3-2x_D^3-4x_Dy_D^2+5x_Dy_D-2x_D)\; x  \\
\quad & + (x_D^4 + x_D^2 y_D^2 - 2 x_D^2 y_D + x_D^2).
\end{aligned}
\end{equation}
\normalsize
The \textbf{factors} command (or the command or \textbf{evala(AFactor)} (for absolute factorization, as we will see in next Section), proves that this generator is irreducible. Actually, we could have conjectured this, as a reducible cubic polynomial can be factorized only with at least one linear factor, and the experimental work with GeoGebra did not reveal any line.

It is now easy to substitute specific coordinates for the point $D$; the resulting cubic curves are identical to those which are obtained by interactive work with GeoGebra.

\begin{enumerate}
\item Take $D(-2,6)$; this point lies inside the parabola. The equation of the corresponding pedal curve is
  \begin{equation}
 \label{eq pedal of parabola wrt internal point}
 x^2 y - 5x^2 + 2x y - 8x + y^3 - 12y^2 + 36y +4=0
 \end{equation}
 Denote by $F(x,y)$ the left-hand side of Equation (\ref{eq pedal of parabola wrt internal point}). The systems of equations $\frac{\partial{F}}{\partial{x}}= \frac{\partial{F}}{\partial{y}}=0$ has one real solution, namely the point $D$, and non real solutions. All of these are irrelevant, thus the curve has no singular point.
\item  Take now $D(4,4)$, a point on the parabola. The equation of the pedal curve is
  \begin{equation}
 \label{eq pedal of parabola wrt point on the parabola}
 x^2 y - 3x^2 - 4x y + 8x + y^3 - 8y^2 + 16y +16=0
 \end{equation}
  Denote by $F(x,y)$ the left-hand side of Equation (\ref{eq pedal of parabola wrt internal point}). The systems of equations $\frac{\partial{F}}{\partial{x}}= \frac{\partial{F}}{\partial{y}}=0$ has 3 real solutions. Two of them correspond to points which do not lie on the curve. The 3rd one is the point $D$ itself, which is a singular point. In order to show that it is a cusp, it is more convenient to use a parametric representation of the curve. We derive such a parametrization by a classical way: as $D$ is a singular point, we consider a line $l_D$ through $D$ with slope $t$. Intersecting this line with the curve $\mathcal{P}$, we find 
  \small
    \begin{equation}
  \label{param pedal of parabola wrt point on the parabola}
  \begin{cases}
    x_1(t) = \frac{4t^3 - 4t^2 - 1}{t(t^2 + 1)}\\
    y_1(t) = -\frac{4t - 3}{t^2 + 1}
  \end{cases}
  \end{equation}
  which is a rational parametrization of $\mathcal{P}$.
  
 Solving the system of equations $(x_1(t)=y_1(t)=4$, we find that the point $D$ corresponds to $t=-1/2$. Now we differentiate $x_1(t)$ and $y_1(t)$ with respect to $t$, and obtain
 \small
    \begin{equation}
  \label{1st derivatives param - point on the parabola}
  \begin{cases}
\frac{d x_1}{d t} = \frac{4 t^{4}+8 t^{3}-t^{2}+1}{t^{2} \left(t^{2}+1\right)^{2}}\\
 \frac{d y_1}{d t} =\frac{4 t^{2}-6 t-4}{\left(t^{2}+1\right)^{2}}
  \end{cases}
  \end{equation}
  \normalsize
 Denote $\overset{\longrightarrow}{V_1}(t)=\left( \frac{d x_1}{d t} ,\frac{d y_1}{d t} \right)$. We have $\overset{\longrightarrow}{V_1}(4)=\overset{\rightarrow}{0}$, which confirms that $D$ is a singular point of the curve.

Now, denote $\overset{\longrightarrow}{V_2}(t)= \frac{d}{dt} \overset{\longrightarrow}{V_1}(t)$. We have:
\footnotesize
\begin{equation*}
\label{2nd derivatives param - point on the parabola}
\overset{\longrightarrow}{V_2}(t)=\left( \frac{-8 t^{6}-24 t^{5}+12 t^{4}+8 t^{3}-6 t^{2}-2}{t^{3} \left(t^{2}+1\right)^{3}}, \frac{-8 t^{3}+18 t^{2}+24 t-6}{\left(t^{2}+1\right)^{3}} \right)
\end{equation*}
\normalsize
and by substitution obtain that $\overset{\longrightarrow}{V_2}\left( -\frac 12 \right)=\frac{32}{5}(2,-1)$.
We differentiate once again:
\footnotesize
\begin{equation*}
\label{3rd derivatives param - point on the parabola}
\begin{split}
 \overset{\longrightarrow}{V_3}(t) & =\frac{d}{dt} \overset{\longrightarrow}{V_2}(t)\\
\quad & = \left( \frac{24 t^{8}+96 t^{7}-84 t^{6}-96 t^{5}+54 t^{4}+24 t^{2}+6}{t^{4} \left(t^{2}+1\right)^{4}},
 \frac{24 \left(t^{4}-3 t^{3}-6 t^{2}+3 t+1\right)}{\left(t^{2}+1\right)^{4}} \right)
\end{split}
\end{equation*}
\normalsize
and obtain $\overset{\longrightarrow}{V_3}\left( -\frac 12 \right)=\frac{384}{25}(7,-1)$. The fact that the vectors $\overset{\longrightarrow}{V_2}\left( -\frac 12 \right)$ and $\overset{\longrightarrow}{V_3}\left( -\frac 12 \right)$ are linearly independent proves that $D$ is a cusp of the pedal curve under study.
\item We take now $D(-6,2)$, a point out of the parabola. A proof that the pedal curve has a crunode (a point of self-intersection) may be not easy, but an efficient method is explained in \cite{curvature} and has been used \cite{DPT}. Here the equation of the pedal curve is
  \begin{equation}
 \label{eq pedal of parabola wrt external point}
x^2 y+y^3-x^2+6 x y-4 y^2+4 y+36=0.
 \end{equation}

Using the intersection of a general line through $D$ (its equation is $y-2=t(x+6)$ with the pedal curve $\mathcal{P}$ , we obtain the following parametrization for $\mathcal{P}$:
\small
\begin{equation}
\label{eq pedal of parabola wrt point out of the parabola}
\begin{cases}
x = -\frac{6 t^3+2 t^2+1}{t \left(t^2+1\right)}\\
y = \frac{6 t+1}{t^2+1}
\end{cases}
\end{equation}
\end{enumerate}

\begin{figure}[htb]
\begin{center}
 \includegraphics[width=12cm]{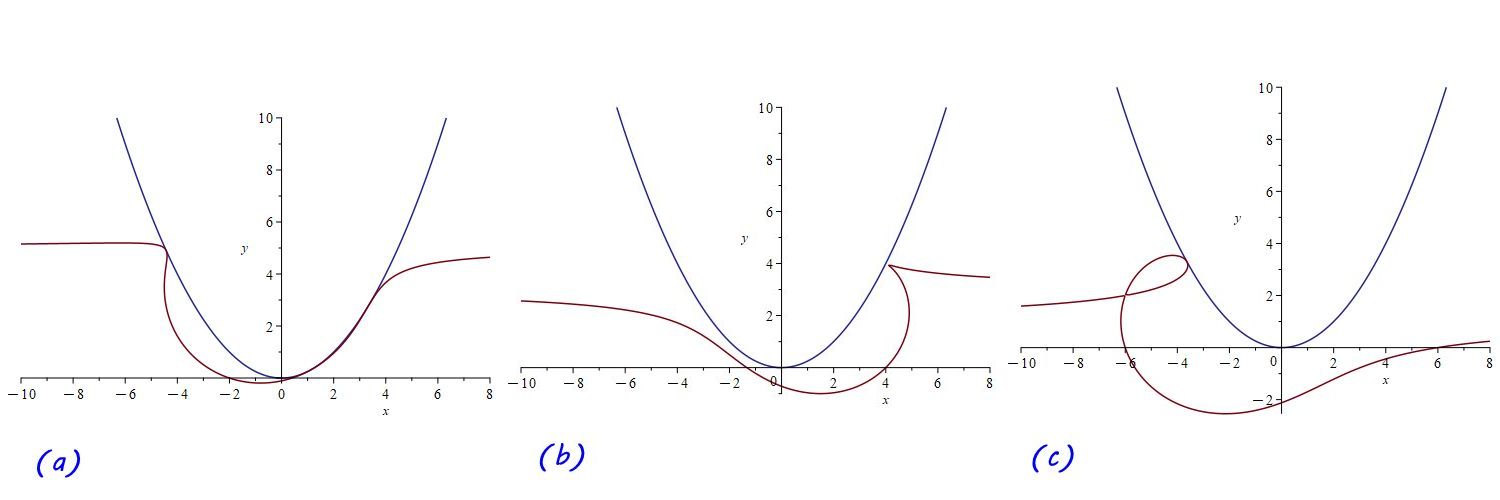}
 \caption{Pedal curve of a canonical parabola w.r.t a point, plotted with Maple}
\label{fig pedal of a parabola - maple}
\end{center}
\end{figure}
Compare the plots in Figures \ref{fig pedal of a parabola - maple} with Figure \ref{fig pedal of a parabola}.
In order to prove that $D$ is a point of self-intersection, we have to find two different values of the parameter corresponding to the same point, i.e., we need to solve the following system of equations:
\begin{equation}
\label{system for self-intersection}
\begin{cases}
  \frac{6 t_1^3+2 t_1^2+1}{t_1 \left(t_1^2+1\right)} = \frac{6 t_2^3+2 t_2^2+1}{t_2 \left(t_2^2+1\right)}\\
  \frac{6 t_1+1}{t_1^2+1} = \frac{6 t_2+1}{t_2^2+1}
\end{cases}
\end{equation}
Maple's \textbf{solve} command yields the following output:
\begin{equation}
(t_1,t_2)=\left( \frac 32 +\frac{\sqrt{7}}{2}, \frac 32-\frac{\sqrt{7}}{2} \right).
\end{equation}
An easy substitution confirms that these values of the parameter correspond both to the point $D$.

The examples above are a good incitement to prove the following theorem:
\begin{theorem}
Let $\mathcal{C}$ be a parabola and $D$ a point in the plane. We denote by $mathcal{P}$ the pedal curve of $mathcal{C}$ with respect to $D$.
\begin{enumerate}
\item If $D$ is a point internal to $\mathcal{C}$, then $\mathcal{P}$ has no singular point.
\item If $D$ is on $\mathcal{C}$, then it belongs also to $\mathcal{P}$ and is a cusp of $\mathcal{P}$.
\item If $D$ is a point external to $\mathcal{C}$, then $\mathcal{P}$ has a point of self-intersection at $D$, and it is the only singular point of $\mathcal{P}$.
\end{enumerate}
\end{theorem}
The proof follows exactly what has been shown in the examples, starting from Equation (\ref{eq pedal of the parabola}). The CAS performs the requested computations, but as the output is really heavy, we do not wish to present it here; the interested reader can run the Maple code.

\section{Pedal curves of a canonical ellipse with respect to an external point}
\label{section ellipses}

We illustrate a general question by means of an example, and display a software-based construction in two different ways.

\subsection{The ellipse is given by a canonical equation}
\begin{figure}[htb]
\begin{center}
 \includegraphics[width=4cm]{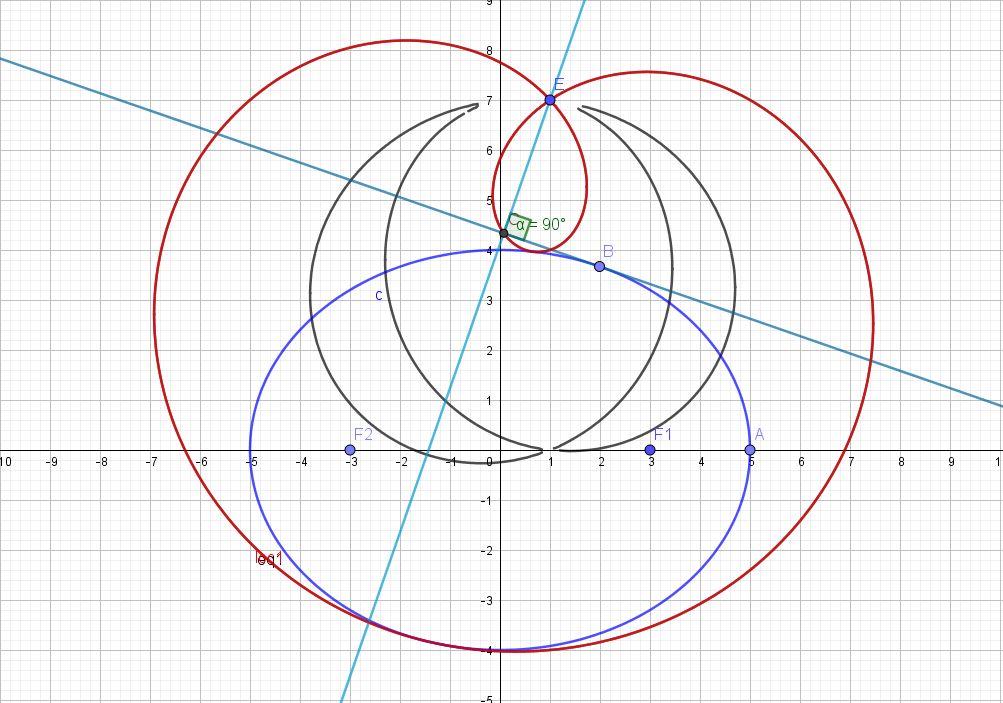}
 \qquad
 \includegraphics[width=5.5cm]{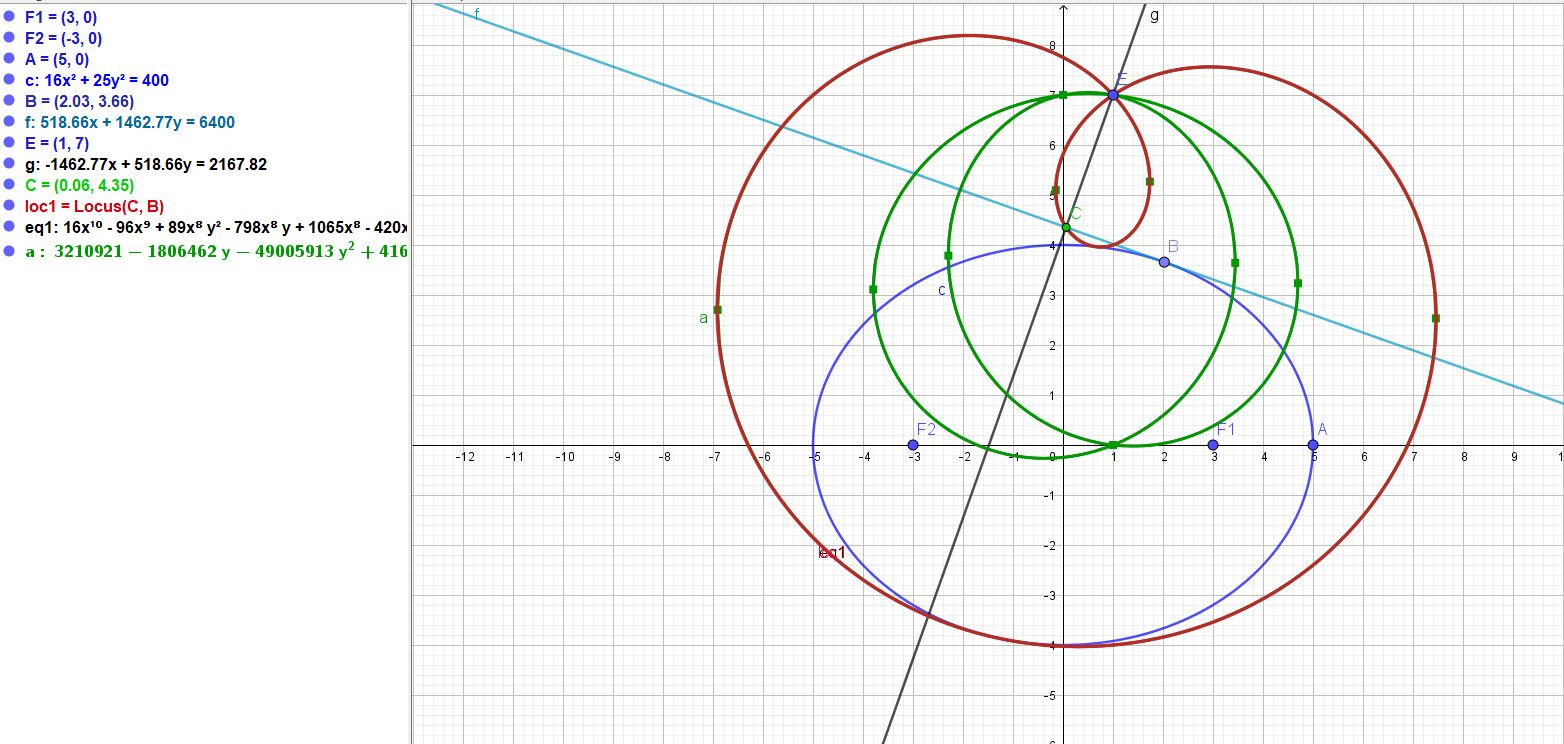}
 \caption{Pedal curve of a canonical ellipse with respect to an external point - total geometric construction}
\label{fig pedal of ellipse 2}
\end{center}
\end{figure}

Let $\mathcal{C}$ be the ellipse whose equation is $\frac{x^2}{25}+\frac{y^2}{16}=1$. Now we follow the following steps:
\begin{enumerate}[(i)]
\item Plot a point $E$ in the plane.
\item Plot the two foci $F1$ and $F2$, and a point $A$ in the plane; for Figure \ref{fig pedal of ellipse 2}, we chose $F1(3,0)$, $F2(-3,0)$ and $A(5,0)$.
\item With GeoGebra command for the building of an ellipse, when the foci and a point are given, construct the ellipse.
\item Plot a point $B$ on $\mathcal{C}$.
\item Plot the tangent to $\mathcal{C}$ at $B$.
\item Plot the perpendicular to this tangent via $E$. The foot of this perpendicular is denoted by $B$.
\item Use the \textbf{Locus(C,B)} command to determine the pedal curve of $\mathcal{C}$ with respect to $E$.
\item Use the \textbf{LocusEquation(C,B)} command to obtain an equation for the pedal curve of $\mathcal{C}$ with respect to $E$.
\end{enumerate}

\begin{figure}[htb]
\begin{center}
 \includegraphics[width=5.5cm]{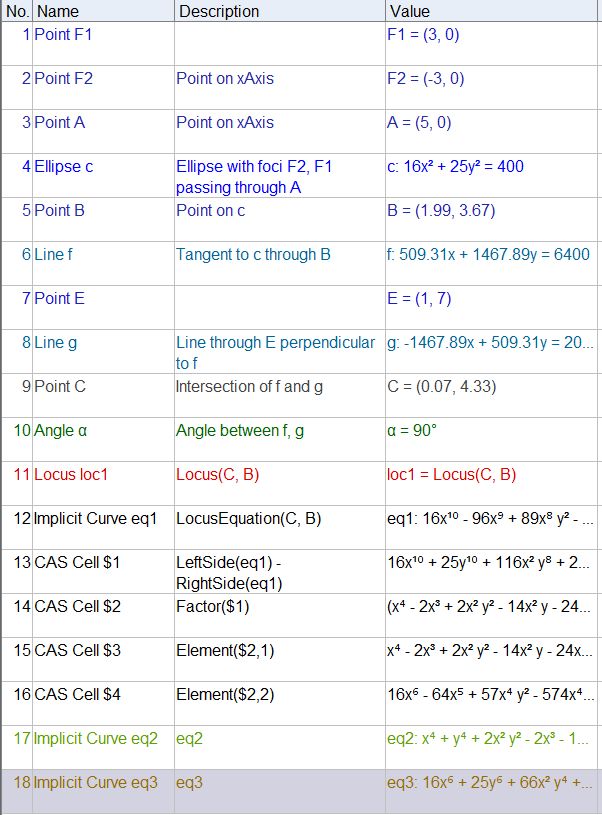}
 \caption{Construction Protocol of the pedal curve}
\label{fig protocol pedal of ellipse 1}
\end{center}
\end{figure}

The graphical output of the last commands are not identical; note that they have been obtained by two different algorithms. In order to have more insight, we switch to the CAS window in order to factorize the polynomial of order 10 which has been displayed. The output of a GeoGebra session is displayed in Figure \ref{fig protocol pedal of ellipse 1}; the Construction Protocol displays the session according to the actual order of commands.

The \textbf{LocusEquation(C,B)} command yields a much more complicated plot and a polynomial equation $P(x,y)=0$, where:
\footnotesize
\begin{equation}
\label{eq pedal of ellipse 1}
\begin{aligned}
P(x,y) = & 16 \; x^{10} + 25 \; y^{10} + 116 \; x^{2} \; y^{8} + 214 \; x^{4} \; y^{6} + 196 \; x^{6} \; y^{4}  \\
& + 89 \; x^{8} \; y^{2} - 96 \; x^{9} - 1050 \; y^{9}  - 132 \; x \; y^{8} - 3948 \; x^{2} \; y^{7} - 492 \; x^{3} \; y^{6}\\
&  - 5544 \; x^{4} \; y^{5} - 684 \; x^{5} \; y^{4} - 3444 \; x^{6} \; y^{3} - 420 \; x^{7} \; y^{2}  - 798 \; x^{8} \; y \\
& + 1065 \; x^{8} + 17991 \; y^{8} + 4396 \; x \; y^{7} + 52284 \; x^{2} \; y^{6} + 12432 \; x^{3} \; y^{5} \\
& + 51660 \; x^{4} \; y^{4}  + 11676 \; x^{5} \; y^{3} + 18432 \; x^{6} \; y^{2} + 3640 \; x^{7} \; y - 2820 \; x^{7} \\
&   - 155148 \; y^{7} - 55748 \; x \; y^{6}  - 318248 \; x^{2} \; y^{5} - 108736 \; x^{3} \; y^{4} - 184912 \; x^{4} \; y^{3} \\
& - 55808 \; x^{5} \; y^{2} - 21812 \; x^{6} \; y - 45548 \; x^{6} + 610117 \; y^{6} + 306376 \; x \; y^{5}\\
&  + 629232 \; x^{2} \; y^{4} + 318920 \; x^{3} \; y^{3} + 1044 \; x^{4} \; y^{2} + 26908 \; x^{5} \; y + 191556 \; x^{5} \\
&  + 227934 \; y^{5}  - 328152 \; x \; y^{4} + 2239888 \; x^{2} \; y^{3} + 411952 \; x^{3} \; y^{2} + 1211700 \; x^{4} \; y  \\
& - 1148762 \; x^{4}  - 11708291 \; y^{4}  - 3671304 \; x \; y^{3} - 12244380 \; x^{2} \; y^{2} - 3039120 \; x^{3} \; y  \\
& + 1714852 \; x^{3} + 41625976 \; y^{3}  + 13815404 \; x \; y^{2}  + 14067452 \; x^{2} \; y + 2699108 \; x^{2}  \\
& - 49005913 \; y^{2} - 10441508 \; x \; y - 6620292 \; x - 1806462 \; y + 3210921.
\end{aligned}
\end{equation}
\normalsize
In the CAS window, we show that $P(x,y)$ is reducible and has two non constant factors $P_1(x,y)$ and $P_2(x,y)$ as follows:
\footnotesize 
\begin{equation}
\label{eq pedal of ellipse 1 - factors}
\begin{aligned}
P_1(x,y) = & \; x^4 - 2   x^3 + 2   x^2   y^2 - 14   x^2   y - 24   x^2 - 2   x  y^2 + 14   x   y + 50   x + y^4 \\
& \; - 14   y^3 + 33   y^2 + 224   y - 809 \\
P_2(x,y) = & \; 16   x^6 - 64   x^6 + 57   x^6   y^6 - 574   x^4   y + 1321   x^4 - 146   x^3   y^2 + 1372   x^3   y \\
& - 2514   x^3 + 66   x^2   y^4 - 1274   x^2   y^3 + 8174   x^2   y^2 - 17038   x^2   y - 2728   x^2 - 82   x   y^4 \\
& + 1498   x   y^3 - 8788   x   y^2 + 15106   x   y + 7938   x + 25   y^6 - 700   y^5 + 7366   y^4  \\
& - 34524   y^3 + 60728   y^2 + 1134   y - 3969.
\end{aligned}
\end{equation}
\normalsize
The equation $P_1(x,y)=0$ determines a quartic, a Lima\c con like curve $\mathcal{P}_1$ (the green curve in Figure \ref{fig pedal of ellipse 2}(b)) and the equation $P_2(x,y)=0$ determines a sextic  $\mathcal{P}_2$ (in red). Actually, the first output is as displayed in Figure \ref{fig pedal of ellipse 2}(a)), it has been corrected with the \textbf{Plot2D} command, available in GD only.
 
What is left to do is to identify that only the quartic  $\mathcal{P}_1$ answers the question. It seems that this can be done only by dragging the point $B$ along the given ellipse, not by algebraic means.

\subsection{Algebraic work with Maple}

As expected, Maple's \textbf{evala(AFactor)} command provides the same factorization as GD. What is left to do is to identify whether the curve $\mathcal{P}_1$ is a Pascal Limac\c con as conjectured, or not. In the previous subsection, we worked on purpose with a point $E$ not on one of the symmetry axes of the ellipse. After all, we worked already with a canonical ellipse. If we had chosen a point $E$ on the $y-$axis, the pedal curve would have been symmetric about the $y-$axis. For example, if $E(0,6)$, the pedal curve has the following quartic equation:

\footnotesize 
\begin{equation}
\label{eq pedal of ellipse 1 - factors}
P_1(x,y) = x^4 + 2 x^2 y^2 - 12 x^2 y - 25x^2 + y^4 - 12  y^3 + 20 y^2 + 192  y- 576
\end{equation}
\normalsize
A canonical equation of a Pascal Lima\c con is given by
\footnotesize 
\begin{equation}
\label{eq of Pascal Limacon}
L(x,y) = (x^2+y^2+a \; e \; y)^2-a^2(x^2+y^2). 
\end{equation}
\normalsize
Let us try to find a change of coordinates of the form
\begin{equation}
\label{change of coordinates}
\begin{cases}
 x= p_1 \; X + q_1 \; Y +r_1\\
 y= p_2 \; X + q_2 \; Y +r_2
\end{cases}.
\end{equation}

The real parameters in Equations (\ref{change of coordinates}) must be determined; we do not need terms in $xy$ as we conjecture that no rotation has to be involved. By substitution into Equation (\ref{eq of Pascal Limacon}), we obtain:
\footnotesize 
\begin{equation}
\label{eq of Pascal Limacon}
\begin{aligned}
L(x,y) = & \left[ (p_1 \; X + q_1 \; Y +r_1)^2+(p_2 \; X + q_2 \; Y +r_2)^2+a \; e \; (p_2 \; X + q_2 \; Y +r_2) \right]^2\\
& -a^2((p_1 \; X + q_1 \; Y +r_1)^2+(p_2 \; X + q_2 \; Y +r_2)^2). 
\end{aligned}
\end{equation}
\normalsize
After expansion, we obtain:
%\tiny
\footnotesize
\begin{equation}
\label{eq of Pascal Limacon 2}
\begin{aligned}
L(x,y) = & (q_1^4 + 2 q_1^2 q_2^2 + q_2^4) Y^4 + (4 p_1 q_1^3 + 4 p_1 q_1 q_2^2 + 4 p_2 q_1^2 q_2 + 4 p_2 q_2^3) XY^3 \\
& + (2 a e q_1^2 q_2 + 4 q_1^3 r_1 + 2 a e q_2^3 + 4 q_1^2 q_2 r_2 + 4 q_1 q_2^2 r_1 + 4 q_2^3 r_2) Y^3 \\
& + (-a^2 q_1^2)Y^2 + ((6 p_1^2 q_1^2 + 2 p_1^2 q_2^2 + 8 p_1 p_2 q_1 q_2 + 2 p_2^2 q_1^2 + 6 p_2^2 q_2^2) X^2Y^2 \\
& + (4 a e p_1 q_1 q_2 + 2 a e p_2 q_1^2 + 6 a e p_2 q_2^2 + 12 p_1 q_1^2 r_1 + 8 p_1 q_1 q_2 r_2 + 4 p_1 q_2^2 r_1 \\
& \qquad + 4 p_2 q_1^2 r_2 + 8 p_2 q_1 q_2 r_1 + 12 p_2 q_2^2 r_2) XY^2 \\
& + (6 a e q_2^2 r_2 + 8 q_1 q_2 r_1 r_2 + 4 a e q_1 q_2 r_1 + 2 a e q_1^2 r_2 + 6 q_1^2 r_1^2 + 2 q_1^2 r_2^2 \\
& \qquad + a^2 e^2 q_2^2 - a^2 q_2^2 + 2 q_2^2 r_1^2 + 6 q_2^2 r_2^2) Y^2\\
& + (4 p_1^3 q_1 + 4 p_1^2 p_2 q_2 + 4 p_1 p_2^2 q_1 + 4 p_2^3 q_2)  X^3Y \\
& + (2 a e p_1^2 q_2 + 4 a e p_1 p_2 q_1 + 6 a e p_2^2 q_2 + 12 p_1^2 q_1 r_1 + 4 p_1^2 q_2 r_2 + 8 p_1 p_2 q_1 r_2\\
& \qquad  + 8 p_1 p_2 q_2 r_1 + 4 p_2^2 q_1 r_1 + 12 p_2^2 q_2 r_2) X^2Y \\
& +  (2 a^2 e^2 p_2 q_2 + 4 a e p_1 q_1 r_2 + 4 a e p_1 q_2 r_1 + 4 a e p_2 q_1 r_1 + 12 a e p_2 q_2 r_2 - 2 a^2 p_1 q_1 \\
& \qquad - 2 a^2 p_2 q_2 + 12 p_1 q_1 r_1^2 + 4 p_1 q_1 r_2^2 + 8 p_1 q_2 r_1 r_2 \\
& \qquad + 8 p_2 q_1 r_1 r_2 + 4 p_2 q_2 r_1^2 + 12 p_2 q_2 r_2^2)  XY \\
& +(- 2 a^2 q_1 r_1 + 2 a^2 e^2 q_2 r_2 + 2 a e q_2 r_1^2 + 6 a e q_2 r_2^2 + 4 q_1 r_1 r_2^2 + 4 a e q_1 r_1 r_2 + 4 q_1 r_1^3 \\
& \qquad - 2 a^2 q_2 r_2 + 4 q_2 r_1^2 r_2 + 4 q_2 r_2^3) Y + (p_1^4 + 2 p_1^2 p_2^2 + p_2^4) X^4\\
& + (2 a e p_1^2 p_2 + 2 a e p_2^3 + 4 p_1^3 r_1 + 4 p_1^2 p_2 r_2 + 4 p_1 p_2^2 r_1 + 4 p_2^3 r_2) X^3 \\
&+ (a^2 e^2 p_2^2 + 2 a e p_1^2 r_2 + 4 a e p_1 p_2 r_1 + 6 a e p_2^2 r_2 - a^2 p_1^2 - a^2 p_2^2 + 6 p_1^2 r_1^2 + 2 p_1^2 r_2^2 \\
& \qquad + 8 p_1 p_2 r_1 r_2 + 2 p_2^2 r_1^2 + 6 p_2^2 r_2^2) X^2 \\
&+ (2 a^2 e^2 p_2 r_2 + 4 a e p_1 r_1 r_2 + 2 a e p_2 r_1^2 + 6 a e p_2 r_2^2 - 2 a^2 p_1 r_1 - 2 a^2 p_2 r_2 + 4 p_1 r_1^3 \\
& \qquad + 4 p_1 r_1 r_2^2 + 4 p_2 r_1^2 r_2 + 4 p_2 r_2^3) X\\
& - a^2 r_1^2 - a^2 r_2^2 + r_2^4 + 2 a e r_2^3 + 2 a e r_1^2 r_2 + 2 r_1^2 r_2^2 + r_1^4 + a^2 e^2 r_2^2.
\end{aligned}
\end{equation}
\normalsize
By identification of the coefficients between Equations (\ref{eq pedal of ellipse 1 - factors}) and (\ref{eq of Pascal Limacon 2}), we have to solve the following system of equations:
\footnotesize
\begin{equation}
\label{system Pascal Limacon}
\begin{cases}
\begin{aligned}
p_1^4 + 2 p_1^2 p_2^2 + p_2^4   &= 1\\
6 p_1^2 q_1^2 + 2 p_1^2 q_2^2 + 8 p_1 p_2 q_1 q_2 + 2 p_2^2 q_1^2 + 6 p_2^2 q_2^2  &=2\\
2 a e p_1^2 q_2 + 4 a e p_1 p_2 q_1 + 6 a e p_2^2 q_2 + 12 p_1^2 q_1 r_1 + 4 p_1^2 q_2 r_2 + 8 p_1 p_2 q_1 r_2\\
+ 8 p_1 p_2 q_2 r_1 + 4 p_2^2 q_1 r_1 + 12 p_2^2 q_2 r_2 & =-12\\
a^2 e^2 p_2^2 + 2 a e p_1^2 r_2 + 4 a e p_1 p_2 r_1 + 6 a e p_2^2 r_2 - a^2 p_1^2 - a^2 p_2^2 + 6 p_1^2 r_1^2 + 2 p_1^2 r_2^2 \\
+ 8 p_1 p_2 r_1 r_2 + 2 p_2^2 r_1^2 + 6 p_2^2 r_2^2 & =-25\\
q_1^4 + 2 q_1^2 q_2^2 + q_2^4  &=0\\
6 a e q_2^2 r_2 + 8 q_1 q_2 r_1 r_2 + 4 a e q_1 q_2 r_1 + 2 a e q_1^2 r_2 + 6 q_1^2 r_1^2 + 2 q_1^2 r_2^2 \\
 + a^2 e^2 q_2^2 - a^2 q_2^2 + 2 q_2^2 r_1^2 + 6 q_2^2 r_2^2 & =20\\
- 2 a^2 q_1 r_1 + 2 a^2 e^2 q_2 r_2 + 2 a e q_2 r_1^2 + 6 a e q_2 r_2^2 + 4 q_1 r_1 r_2^2 + 4 a e q_1 r_1 r_2 + 4 q_1 r_1^3 \\
 - 2 a^2 q_2 r_2 + 4 q_2 r_1^2 r_2 + 4 q_2 r_2^3 & =192\\
 - a^2 r_1^2 - a^2 r_2^2 + r_2^4 + 2 a e r_2^3 + 2 a e r_1^2 r_2 + 2 r_1^2 r_2^2 + r_1^4 + a^2 e^2 r_2^2 & =576
\end{aligned}
\end{cases}
\end{equation}
\normalsize
all the other coefficients in Equation (\ref{eq of Pascal Limacon 2}) being equal to 0. Note that an expanded form of $L(x,y)$ is a huge expression, but to determine the coefficients we do not need the entire expression. 

\begin{figure}[htb]
\begin{center}
 \includegraphics[width=6cm]{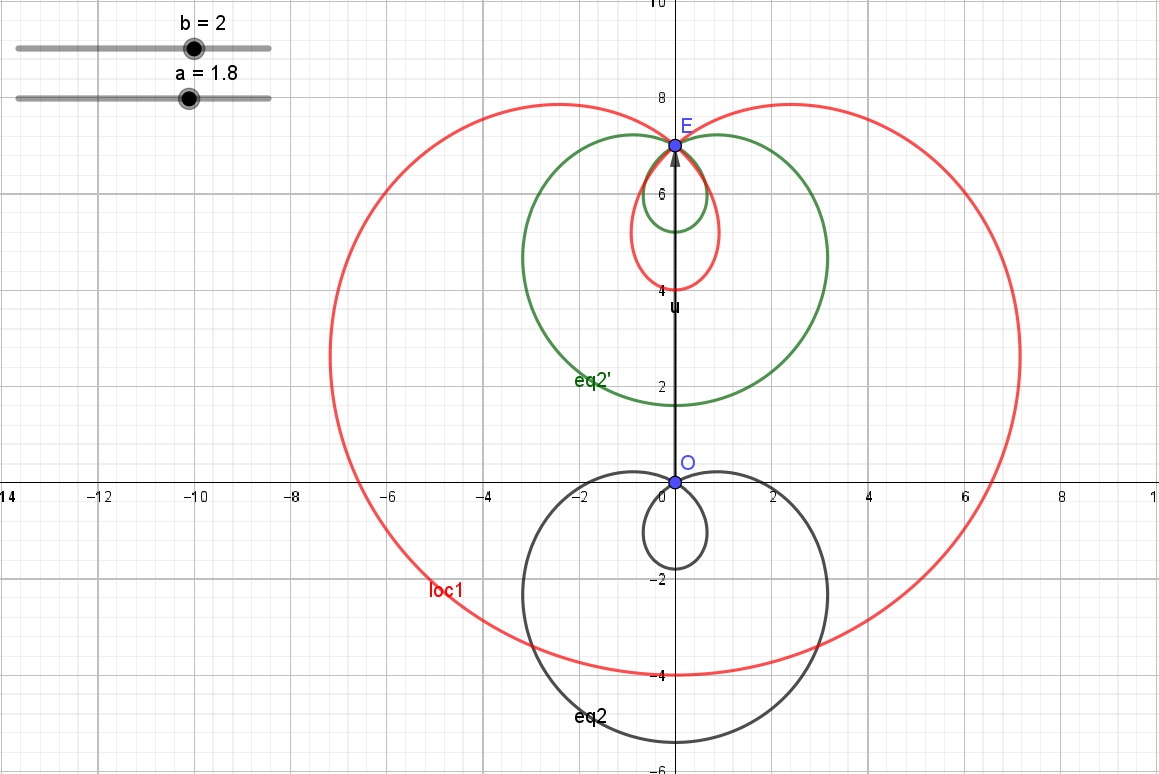}
 \caption{Checking whether we obtained a Lima\c con or not}
\label{fig limacon or not?}
\end{center}
\end{figure}

As the algebraic work is heavy, the explorer may try to have a geometric experimentation, illustrated in Figure \ref{fig limacon or not?}, performed interactively using GD. For this exploration, we chose a point $E$ on the $y-$axis, which ensures that the entire situation is symmetric about the $y-$axis, and enables to use a "canonical"  equation for the lima\c con available in the literature.  Plot the lima\c con given by the equation\footnote{This equation has been adapted from the one given in \url{https://mathcurve.com/courbes2d.gb/limacon/limacon.shtml}.}
\begin{equation}
\label{general eq limacon}
(x^2+y^2+e*a*y)^2-a^2*(x^2+y^2)=0.
\end{equation}
As the letter e has a specific meaning in GD, we used b instead.
Now denote $\overset{\rightarrow}{u}=\overset{\longrightarrow}{OE}$ and translate the plotted lima\c con by this vector. Now change the values of the parameters with the sliders in order to try to have the translated lima\c con (in black) coalesce with the pedal curve obtained previously (in red). Of course, not obtained a true coalescing is not a complete proof, but may be enough for the explorer. Such a situation has been already met in \cite{hyperbolisms}, where curves "looking like but different from" an Agnesi witch have been obtained. There, observing the equations was an important part of the conviction (not of the proof).

\section{Pedal curves of a canonical ellipse with respect to the origin}
\begin{figure}[htb]
\begin{center}
 \includegraphics[width=5cm]{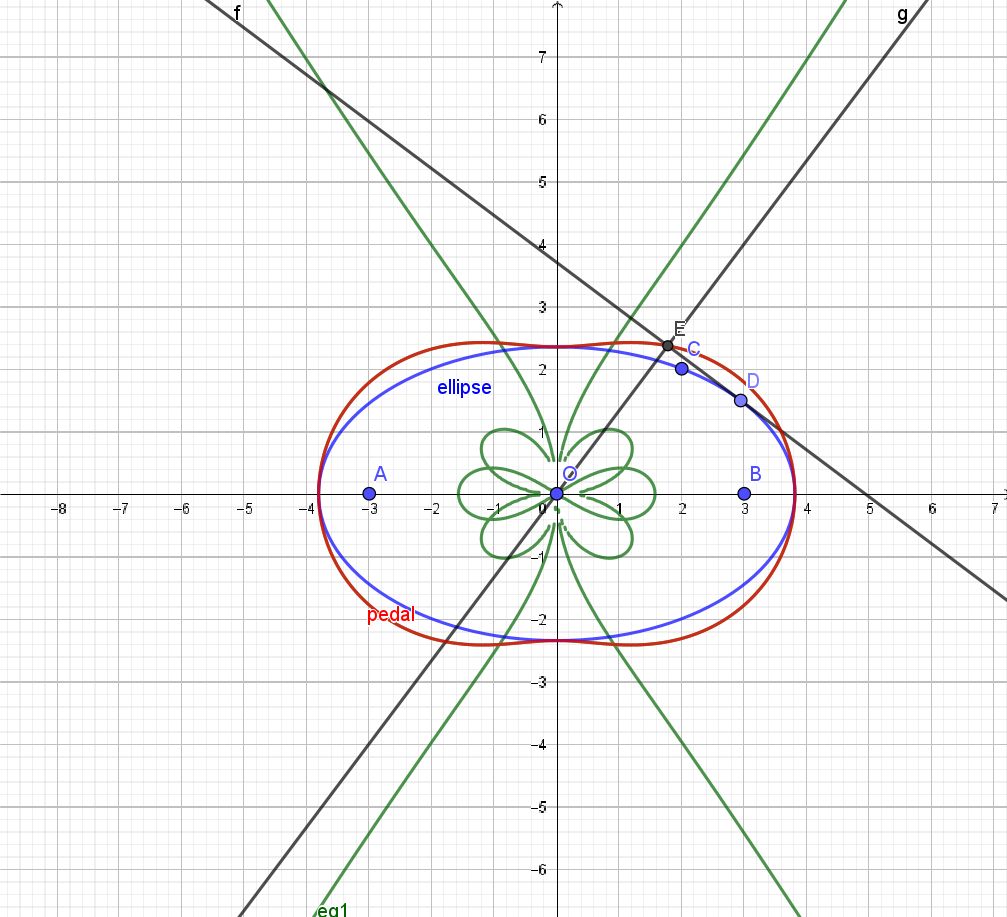}
 \qquad
 \includegraphics[width=5cm]{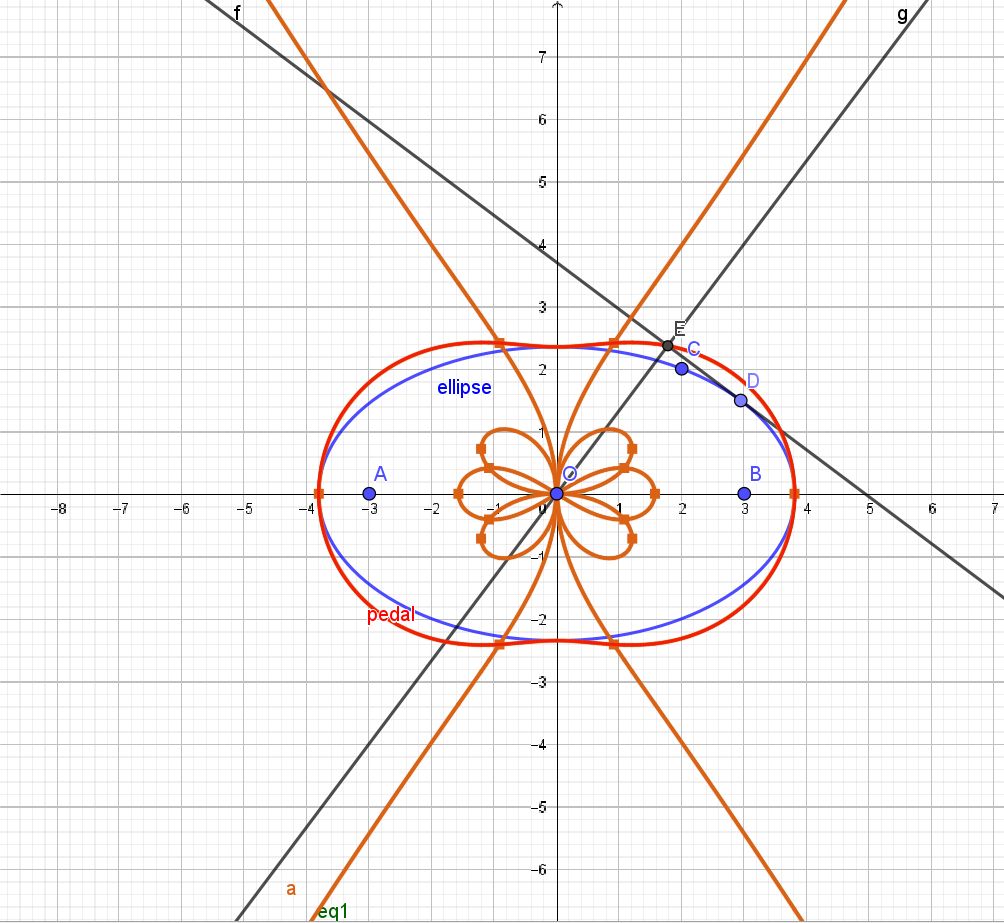}
 \caption{Pedal curve of a canonical ellipse}
\label{fig pedal of an ellipse}
\end{center}
\end{figure}

The \textbf{LocusEquation} command of GD provides the following polynomial, defining the curve plotted in Figure \ref{fig pedal of an ellipse}(a):
\scriptsize
\begin{equation}
\label{eq pedal of ellipse wrt origin}
\begin{aligned}
f(x,y)  = & \quad 4x^{20} + 41x^{18}y^2 + 180x^{16}y^4 + 444x^{14}y^6 + 672x^{12}y^8 + 630x^{10}y^{10}  + 336x^8y^{12} \\
\qquad &   + 60x^6y^{14} - 36x^4y^{16} - 23x^2y^{18}  - 4y^{20} - 68x^{18} - 566x^{16}y^2 - 1738x^{14}y^4 \\ 
 \qquad &   - 2354x^{12}y^6 - 754x^{10}y^8    + 1762x^8y^{10} + 2354x^6y^{12} + 1162x^4y^{14} + 206x^2y^{16} \\
\qquad &    - 4y^{18} + 144x^{16} 729x^{14}y^2     - 2880x^{12}y^4 - 15777x^{10}y^6 - 24624x^8y^8 - 15777x^6y^{10}\\
\qquad &   - 2880x^4y^{12} + 729x^2y^{14} + 144y^{16} - 324x^{12}y^2 + 17982x^{10}y^4 + 24138x^8y^6  \\
\qquad &   - 12474x^6y^8  - 23814x^4y^{10} - 5508x^2y^{12} - 26244x^8y^4 + 59049x^6y^6 + 26244x^4y^8\\
\end{aligned}
\end{equation}
\normalsize
Switching to the CAS window of GD, we apply the \textbf{factor} command. The output of the \textbf{factor} command is a product of two polynomials, one of degree 12 and one of degree 8, namely:
\scriptsize
\begin{align*}
P_1(x,y)= & 4x^{12} + 25x^{10}y^2 + 56x^8y^4 + 54x^6y^6 + 16x^4y^8 - 7x^2y^{10}  \\
\quad & - 4y^{12} - 9x^8y^2 + 135x^6y^4 + 297x^4y^6 + 153x^2y^8 - 729x^4y^4 \\
P_2(x,y)= & x^{8}+4 x^{6} y^{2}+6 x^{4} y^{4}+4 x^{2} y^{6}+y^{8}-17 x^{6}\\
\quad & -33 x^{4} y^{2}-15 x^{2} y^{4}+y^{6}+36 x^{4}-81 x^{2} y^{2}-36 y^{4} 
\end{align*}
\normalsize
This means already that the obtained curve is reducible, but as we will see this decomposition is not the ultimate one. 

\begin{figure}[htb]
\begin{center}
 \includegraphics[width=7cm]{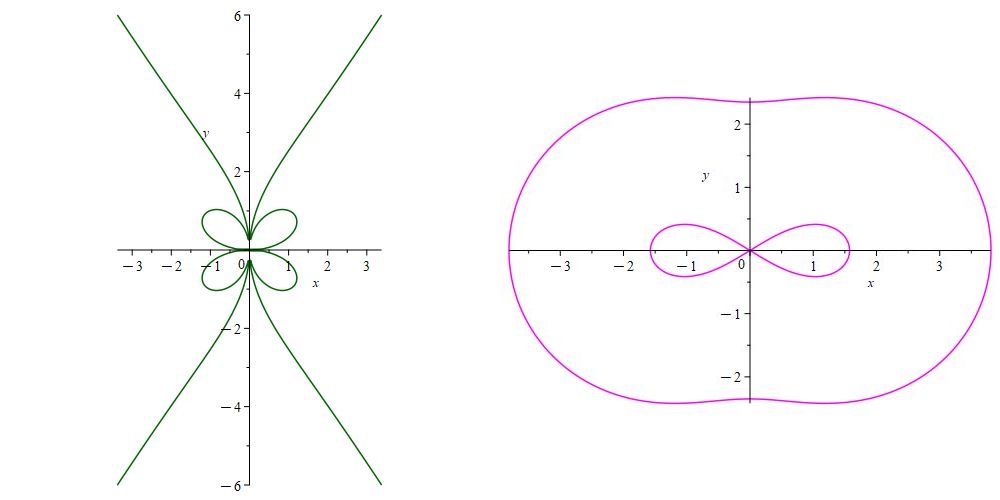}
 \caption{Separate plots with Maple}
\label{fig pedal of an ellipse - maple}
\end{center}
\end{figure}

We switch now to Maple; its \textbf{factor} command provides the same answer as Giac. Separate plots for these polynomials are displayed 
in Figure \ref{fig pedal of an ellipse - maple}.
But we could go further, with a more advanced Maple command for factorization (as in \cite{octics}).
 \small

 \begin{verbatim}
evala(AFactor(f));
F := allvalues(%);
 \end{verbatim}
 \normalsize
 
The first command provides an answer involving the place holder \emph{RootOf}, which can be resolved using the \textbf{allvalues} command. The output provides the following polynomial, whose components are of multiplicity 2.
\scriptsize
\begin{align*}
F=& 4 \; \left[ x^{4}+\left(2 y^{2}-\frac{\sqrt{145}}{2}-\frac{17}{2}\right) x^{2}+y^{4}+\left(-\frac{\sqrt{145}}{2}+\frac{1}{2}\right) y^{2} \right]^2 \\
\quad &
\left[ x^{4}+\left(2 y^{2}+\frac{\sqrt{145}}{2}-\frac{17}{2}\right) x^{2}+y^{4}+\left(\frac{\sqrt{145}}{2}+\frac{1}{2}\right) y^{2} \right]^2\\
\quad &
\left[ x^{6}+\left(-\frac{\sqrt{145}}{8}+\frac{25}{8}\right) y^{2} x^{4}+\left(\left(-\frac{\sqrt{145}}{4}+\frac{13}{4}\right) y^{4}+\left(\frac{9 \sqrt{145}}{8}-\frac{9}{8}\right) y^{2}\right) x^{2}+\left(-\frac{\sqrt{145}}{8}+\frac{9}{8}\right) y^{6} \right]^2 \\
\quad &
\left[ x^{6}+\left(\frac{\sqrt{145}}{8}+\frac{25}{8}\right) y^{2} x^{4}+\left(\left(\frac{\sqrt{145}}{4}+\frac{13}{4}\right) y^{4}+\left(-\frac{9 \sqrt{145}}{8}-\frac{9}{8}\right) y^{2}\right) x^{2}+\left(\frac{\sqrt{145}}{8}+\frac{9}{8}\right) y^{6} \right]^2
\end{align*}
\normalsize
Here too, an interactive exploration may provide conviction which part of the plot is relevant to the actual question. Once again, irrelevant components appear for topological reasons, as in \cite{new approach}.

\section{Discussion}
Pedal curves are a classical topic in plane geometry. Many catalogues are devoted to them \cite{fischer,Lawrence catalog,yates,salmon}. They are part of the vast topic provided by the search and study of geometric loci; we refer here to small sample of papers such as to  \cite{blazek-pech,thales etc.}. In this work, we could shed a new light on the topic. One of the main products is the discovery and study of algebraic curves of higher degree than what is traditionally offered to students. In the recent past, such studies yielded new contributions to the study of quartics (in particular Cassini ovals in \cite{DMZ} and others \cite{DPK-2018,isoptics of quartics}) and the study of sextics \cite{DPT}and octics \cite{octics}. 

Generally, working in a technology-rich environment enables to proceed according a non-traditional frame, namely exploration-conjecture-proof. Slowly, this becomes a main stream way of working in geometry. Sendra and Winkler's book \cite{sendra and winkler} is devoted to parametric curves with a Computer Algebra approach. Automated methods are under constant development \cite{recio et al}, as we saw with GeoGebra-Discovery, and other packages are available. Once again here, as in \cite{dialog,thales etc. MT} we face a situation where a single package of software does not provides the whole picture that the user is trying to build. in the present state-of-the-art, the data transfer from one software to the other is made "by hand" (copy-paste and then fine tuning to adapt what is written to the specific syntax of the receiving software). Following Roanes-Lozano's wish \cite{eugenio}, we hope that this dialog will become more and more automatic.
 
Of course, the software in use depends strongly on the needed affordances, but also on the institutional culture and decisions of the policy makers (see \cite{artigue}). The work presented in this paper has been proposed to in-service teachers learning towards an advanced degree in Mathematics Education, together with topics such as \cite{revival,hyperbolisms}. Their feedback was enthusiastic, as they discovered that mathematics is a very dynamic field.  

% ------------------------------------------------------------------------

\subsection*{Acknowledgment}
The author acknowledge partial support by the CEMJ Chair at JCT.

% ------------------------------------------------------------------------
\end{document}